\documentclass{article}
\usepackage{graphicx,xcolor}
\usepackage{amssymb}
\usepackage[geometry]{ifsym}
\usepackage{rotating}

\newtheorem{lemma}{Lemma}
\newtheorem{theorem}[lemma]{Theorem}

\newcommand{\DONOTTEX}[1]{}

\begin{document}

\title{{\bf Rooted complete minors in \\ line graphs with a Kempe coloring}}
\author{Matthias Kriesell \and Samuel Mohr\footnote%
{Gef\"ordert durch die Deutsche Forschungsgemeinschaft (DFG) -- 327533333.}
}

\maketitle

\begin{abstract}
  \setlength{\parindent}{0em}
  \setlength{\parskip}{1.5ex}

  It has been conjectured that if a finite graph has a vertex coloring such that the union of any two color classes induces a connected graph,
  then for every set $T$ of vertices containing exactly one member from each color class there exists a complete minor such that $T$ contains
  exactly one member from each branching set. Here we prove the statement for line graphs.
  

  
\end{abstract}

\maketitle

Hadwiger’s Conjecture states that the order $h(G)$ of a largest clique minor in a graph $G$ is at least its chromatic number $\chi(G)$~\cite{Hadwiger1943}.
It is known to be true for graphs with $\chi(G) \leq 6$, where for $\chi(G)=5$ or $\chi(G)=6$, we have equivalence to the Four-Color-Theorem, respectively~\cite{RobertsonSeymourThomas1993}.
Instead of restricting the {\em number of color classes}, one could also uniformly bound the {\em order of the color classes}, but even when forbidding anticliques of order $3$ (which bounds
these orders by $2$), the problem is wide open (cf.~\cite{Seymour2016}). In~\cite{Kriesell2016}, the first author suggested to bound the {\em number of colorings}, in particular to consider
uniquely optimally colorable graphs; if $k$ is uniquely $k$-colorable and $x_1,\dots,x_k$ have different colors, then it is easy to see that there exists
a system of edge-disjoint $x_i$,$x_j$-paths ($i \not= j$ from $\{1,\dots,k\}$), a so-called {\em (weak) clique immersion} of order $k$ at $x_1,\dots,x_k$, and the question suggests
itself whether there exists a {\em clique minor} of the same order such that $x_1,\dots,x_k$ are in different bags.
This has been answered affirmatively in~\cite{Kriesell2016} if one forbids antitriangles in $G$. The present paper gives an affirmative answer in the case that $G$ is a line graph.
It should be mentioned that Hadwiger’s Conjecture is known to be true for line graphs in general by a result of Reed and Seymour~\cite{ReedSeymour2004},
but it seems that their argument leaves no freedom for prescribing vertices in the clique minor at the expense of forcing any pair of color classes to be connected.

\medskip

\centerline{*}

All {\em graphs} considered here are supposed to be finite, undirected, and loopless. They may contain parallel edges,
graphs without these are called {\em simple}. For graph terminology not defined here, we refer to~\cite{Diestel2017}. 
A {\em clique minor} of $G$ is a set of connected, nonempty, pairwise disjoint, pairwise adjacent subsets of $V(G)$ (the {\em branching sets}),
where a set $A \subseteq V(G)$ is {\em connected} if the subgraph $G[A]$ induced by $A$ in $G$ is connected,
and disjoint $A,B \subseteq V(G)$ are {\em adjacent} if some vertex of $A$ is adjacent to some vertex of $B$.
An {\em anticlique} of $G$ is a set of pairwise nonadjacent vertices.
A {\em coloring} of a graph $G$ is a partition $\mathfrak{C}$ into anticliques, the {\em color classes},
and we call it a {\em Kempe coloring} if the union of any two of these induces a connected subgraph in $G$.
Throughout, a (minimal) {\em transversal} of a set $\mathfrak{S}$ of pairwise disjoint sets is a set $T \subseteq \bigcup \mathfrak{S}$ 
such that $|T \cap A|=1$ for all $A \in \mathfrak{S}$; in this case, we also say that $\mathfrak{S}$ is {\em traversed by} $T$.
The {\em line graph} $L(H)$ of some graph $H$ is the (simple) graph with vertex set $E(H)$ where two distinct vertices are adjacent if and only
they are incident (as edges) in $H$.

The first author conjectured in~\cite{Kriesell2016} that
for every transversal $T$ of every Kempe coloring of a graph $G$ there exists a complete minor in $G$ traversed by $T$.
Here we prove the conjecture for line graphs.

\begin{theorem}
  \label{T1}
  For every transversal of every Kempe coloring of the line graph $L(H)$ of any graph $H$ there exists a complete minor
  in $L(H)$ traversed by $T$.
\end{theorem}

Of course this statement can be fomulated entirely without addressing to line graphs. Call a set $F$ of edges of a graph $H$
{\em connected} if any two of them are on a path of edges from $F$, and call two sets $F,F'$ of edges {\em incident} if
some edge of $F$ is incident with some edge of $F'$. Theorem~\ref{T1} translates as follows:

\begin{theorem}
  \label{T2}
  Let $H$ be a graph and $\mathfrak{C}$ be a partition of $E(H)$ into (not necessarily maximum) matchings such that
  the union of any two of them is connected. Then for every transversal $T$ of $\mathfrak{C}$ there
  exists a set of connected, pairwise disjoint, pairwise incident edge sets traversed by $T$.
\end{theorem}

\medskip

\centerline{*}

As our proof of Theorem~\ref{T2} uses contraction at some places, it is reasonable to allow multiple edges (but no loops) in $H$.
However, the precondition of Theorem~\ref{T2} imposes a very special structure on $H$ as soon as $H$ contains a pair of parallel
edges. We thus prefer to give a separate, simple proof for this situation instead of handling parallel edges in the proof of Theorem~\ref{T2}.
Given a graph $H$, let us say that $F \subseteq E(H)$ {\em covers} $v \in V(H)$ if $v$ is incident with at least one edge from $F$.
By $E_H(v)$ we denote the set of all edges incident with $v \in V(H)$.

\begin{lemma}
  \label{T3}
  Let $H$ be a graph with a pair of parallel edges and $\mathfrak{C}$ be a partition of $E(H)$ into (not necessarily maximum) matchings such that
  the union of any two of them is connected. Then for every transversal $T$ of $\mathfrak{C}$, there
  exists a set of connected, pairwise disjoint, pairwise incident edge sets traversed by $T$.
\end{lemma}

{\bf Proof.}
Let $e,f$ be parallel edges of $H$.
They are in different matchings $M_e,M_f$ of $\mathfrak{C}$, and they form a cycle of length $2$,
so that $M_e=\{e\}$ and $M_f=\{f\}$.  Let $M \in \mathfrak{C} \setminus \{M_e,M_f\}$. As every edge in $M$ is incident with some
edge of $M_e$, every edge from $M$ (and hence every edge from $H$) is incident with $e,f$, implying that $M$ contains at most two edges.
If $|M|=1$, say, $M=\{g\}$, then, likewise, $g$ is incident with every other edge of $H$ and
we apply induction to $H-g$, $\mathfrak{C} \setminus \{M\}$, and $T \setminus \{g\}$ as to find
a set $\mathfrak{K}$ of connected, pairwise disjoint, pairwise incident edge sets traversed by $T$,
and $\mathfrak{K} \cup \{\{g\}\}$ proves the statement for $H$ (and $\mathfrak{C},T$). ---
So we may assume that every $M$ distinct from $M_e,M_f$ consists of two edges. In particular, $e,f$ is the only pair of parallel edges in $H$.
Let $x,y$ be the endvertices of $e,f$, and let $M_i=\{xa_i,yb_i\}$, $i \in \{1,\dots,\ell\}$ be the matchings from $\mathfrak{C} \setminus \{M_e,M_f\}$.
The $a_i$ are pairwise distinct, and so are the $b_i$. For $i \not= j$ from $\{1,\dots,\ell\}$,
at least one of $a_i=b_j$ or $b_i=a_j$ holds since $M_i \cup M_j$ is connected.
Since $a_i,b_i$ have degree at most $2$, for $i \in \{1,\dots,\ell\}$,
there exists at most one $j \not= i$ with $a_i=b_j$ and at most one $j \not= i$ with $b_i=a_j$.
It follows that $\ell \leq 3$. We may assume that the $\ell+2$ transversal edges from $T$ are not pairwise incident,
as in this case $\{\{g\}:\,g \in T\}$ proves the statement.
In particular, $\ell \geq 2$, and for $i \not= j$ only one of $a_i=b_j$ and $b_i=a_j$ holds, as otherwise 
$M_i \cup M_j$ induce a $4$-cycle, $\ell=2$, and the four transversal edges are pairwise incident.
Hence, we may assume that if $\ell=2$, then $b_1=a_2$ and $T=\{e,f,xa_1,yb_2\}$ without loss of generality,
and if $\ell=3$, then $b_1=a_2,b_2=a_3,b_3=a_1$ and $T=\{e,f,xa_1,yb_2,yb_3\}$ without loss of generality.
In either case, $\{\{e\},\{f\},\{xa_1,xa_2,yb_1\},\{yb_2\},\dots,\{yb_\ell\}\}$ proves the statement.
\hspace*{\fill}$\Box$

If $\mathfrak{C}$ is the only coloring of order $k$ in a graph $G$, then $\mathfrak{C}$ is a Kempe coloring,
and there is no coloring with fewer than $k$ colors (so $k$ equals the chromatic number $\chi(G)$ of $G$).
Even in this very special case, we do not know whether $G$ has a complete minor of order $k$ (disregarding transversals),
that is, we do not know if Hadwiger's conjecture~\cite{Hadwiger1943} is true for uniquely $k$-colorable graphs.
However, the situation for uniquely $k$-colorable line graphs is almost completely trivial, as for $k \geq 4$, the star $K_{1,k}$
is the only uniquely $k$-edge-colorable simple graph~\cite{Thomason1978} (and the non-simple ones are covered by Lemma~\ref{T3}).

At this point, one may wonder if the graphs considered in Theorem~\ref{T2} are ``rare''. Let us show that this is not the case.
A partition of the edge set of a graph into (perfect) matchings such that the union of any two of them induces a Hamilton cycle in $G$
is called a {\em perfect $1$-factorization}. Clearly, every graph with a perfect $1$-factorization is $k$-regular and 
meets the assumptions of Theorem~\ref{T2}.
There is an old conjecture by Kotzig stating that every complete graph of even order has a perfect $1$-factorization
\cite{Kotzig1963}, indicating that it is difficult to determine whether a graph has a perfect $1$-factorization.
However, for our purposes it suffices to construct some variety of graphs which have one, as done in (A),(B) below.

(A) Given $k$, let $a=(a_1,\dots,a_k)$ be a sequence of
pairwise distinct nonnegative integers and take $m$ larger than all of these and relatively prime to any possible difference
$a_i-a_j$, $i \not= j$; for example, just take a large prime number. Let $V=\mathbb{Z}_m \times \{0,1\}$,
set $M_i:=\{(z,0)(z+\overline{a_i},1):\,z \in \mathbb{Z}_m\}$, where $\overline{a}$ denotes the residual class modulo $m$ containing $a \in \mathbb{Z}$,
and let $E:=M_1 \cup \dots \cup M_k$. The resulting graph $(V,E)=:H(m,a)=:H$ is bipartite, $k$-regular, and $\{M_1,\dots,M_k\}$ is a partition of $E$ into
perfect matchings. The graph induced by the union of $M_i \cup M_j$ thus decomposes into cycles. Take a vertex $(z,0)$ from such a cycle $C$;
by following the $M_i$-edge, we reach $(z+\overline{a_i},1)$, and by following the $M_j$-edge from there we
reach $(z+\overline{a_i}-\overline{a_j},0)$. So $C$ contains all vertices of the form $(z+\overline{t \cdot (a_i-a_j)},0)$, $t \in \mathbb{Z}$,
from $\mathbb{Z}_m \times \{0\}$, and since $a_i-a_j$ and $m$ are relatively prime, it contains {\em all} vertices from $\mathbb{Z}_m \times \{0\}$.
Consequently, $C$ is a Hamilton cycle, so $H$ admits a perfect $1$-factorization of order $k$.

(B) For $i \in \{1,2\}$, take any graph $H_i$ with a perfect $1$-factorization $\mathfrak{C}_i$ of order $k$.
Take a vertex $v_i$ in $V(H_i)$ and consider the edges $e_{i,1},\dots,e_{i,k}$ incident
with $v_i$. Let $M_{i,j}$ be the member from $\mathfrak{C}_i$ containing $e_{i,j}$. Assuming that $H_1,H_2$ are disjoint,
let $H$ be obtained from the union of $H_1-v_1$ and $H_2-v_2$ by adding a new edge $f_j$ from the endvertex of $e_{1,j}$ distinct from $v_1$
to the endvertex of $e_{2,j}$ distinct from $v_2$ for each $j \in \{1,\dots,k\}$.
Then $M_j:=((M_{1,j} \cup M_{2,j}) \setminus \{e_{1j},e_{2,j}\}) \cup \{f_j\}$ for $j \in \{1,\dots,k\}$
defines a perfect matching of $H$, and the union of any two of these induces a Hamilton cycle of $H$,
so $H$ has a perfect $1$-factorization.

(C) If one deletes a vertex of any graph with a perfect $1$-factorization,
then the edge set of the resulting graph has an (obvious) partition into matchings such that the union of any two of them is a Hamilton path;
so we get further graphs meeting the assumptions of Theorem~\ref{T2} this way.

\bigskip 

\centerline{*}

Back to Theorem~\ref{T2},
let us now consider the case that $H$ is a complete graph. It turns out that for {\em any} set $T$ of $n$ edges
(not necessarily being a transversal of some set of matchings as in Theorem~\ref{T2})
we can find connected, pairwise disjoint, pairwise incident edge sets traversed by $T$.

\begin{lemma}
  \label{T4}
  For every set $T$ of $n$ edges of the simple complete graph $H$ on $n \geq 3$ vertices, there exists a set of connected, pairwise disjoint,
  pairwise incident edge sets traversed by $T$.
\end{lemma}

{\bf Proof.}
For $n=3$, the statement is obviously true. For $n>3$, consider the subgraph $H[T]:=(V(H),T)$ induced by $T$.
It has average degree $2$ and, therefore, a vertex $v$ with at most two neighbors in $H[T]$.

If $v$ has exactly one neighbor $x$ in $H[T]$, then we apply induction to $H-v$ and find a set $\mathfrak{K}$ of $n-1$ connected, pairwise disjoint,
pairwise incident edge sets traversed by $T \setminus \{vx\}$, and $\mathfrak{K} \cup \{E_H(v)\}$ proves the statement for $H$.

If $v$ has exactly two neighbors $x,y$ in $H[T]$, then we may assume that $x$ is not incident with all the $n-2$ edges from $T \setminus \{vx,vy\}$,
since otherwise these would form a spanning star in $H-v$ and one of the leaves of this star had degree $1$ in $H[T]$ since $n>3$ --- a case
which we have just considered. Therefore, there exists an edge $xz$ in $E(H-v) \setminus T$.
Induction applied to $H-v$ provides a set $\mathfrak{K}$ of $n-1$ connected, pairwise disjoint, pairwise incident edges traversed
by $(T \setminus \{vx,vy\}) \cup \{xz\}$. Let $F$ be the member of $\mathfrak{K}$ containing $xz$.
The set $E_H(v) \setminus \{vx\}$ is incident to all of $\mathfrak{K}$ in $H$, as each of these cover at least two neighbors of $v$,
so that $(\mathfrak{K} \setminus \{F\}) \cup \{F \cup \{vx\},E_H(v) \setminus \{vx\}\}$ proves the statement for $H$.

If $v$ has no neighbors in $H[T]$ and $xy$ is any edge in $T$, then by induction there exists a set $\mathfrak{K}$ of
$n-1$ connected, pairwise incident edges traversed by $T \setminus \{xy\}$ in $H-v$.
If $xy$ is not an edge of any member of $\mathfrak{K}$,
then $\mathfrak{K} \cup \{E_H(v) \cup \{xy\}\}$ proves the statement for $H$. Otherwise, there is an $F \in \mathfrak{K}$ with $xy \in F$.
Let $wz \not=xy$ be the edge from $T$ contained in $F$, where we may assume that $w \not\in \{x,y\}$.
By symmetry, we may assume that $w$ is in the component of $H[F]-xy$ containing $x$, so that $F':=(F \setminus \{xy\}) \cup \{vw,vy\}$ is
connected and covers $v$ and all vertices covered by $F$. The set $F'':=(E_H(v) \setminus \{vw,vy\}) \cup \{xy\}$ is connected and covers all
vertices of $H$ except for $w$, so that it is incident with $F'$ and all sets from $\mathfrak{K}$ as each of these cover at least two neighbors of $v$.
Consequently, $(\mathfrak{K} \setminus \{F\}) \cup \{F',F''\}$ proves the statement for $H$.
\hspace*{\fill}$\Box$

It remains open if Lemma~\ref{T4} is best possible in the sense that we cannot prescribe a set $T$ of more than $n$ edges there.
For $n=3$ and $n=4$, optimality is easy to check, for $n=5$, one cannot prescribe six edges if they form a subgraph $K_{2,3}$.
In general, for $n>2$, one cannot prescribe $2n-2$ edges if they form a graph $K_{2,n-2}$ plus two edges, one of them connecting the
two vertices of degree $n-2$ in $K_{2,n-2}$, but there should be better bounds.

\bigskip

\centerline{*}

Finishing preparation for the proof of our main result, let us recall and slightly extend Lemmas from~\cite{Kriesell2016} and~\cite{Kriesell2001}.

\begin{lemma}~\cite{Kriesell2016}
  \label{L1}
  Suppose that $\mathfrak{C}$ is a Kempe coloring of order $k$ of a graph $G$ and let $S \subseteq V(G)$ be a separating set. Then
  (i) if $F \in \mathfrak{C}$ does not contain any vertex from $S$, then it contains a vertex from every component of $G-S$, and
  (ii) $S$ contains vertices from at least $k-1$ members of $\mathfrak{C}$.
\end{lemma}

{\bf Proof.}
Let $C,D$ be distinct components of $G-S$ and let $F \in \mathfrak{C}$ with $F \cap S= \emptyset$.
Suppose, to the contrary, that  $F \cap V(C) = \emptyset$. 
Take any vertex $x \in V(C)$ and the set $A$ from $\mathfrak{C}$ containing $x$.
Then $x$ has no neighbors in $F$, so that $G[A \cup F]$ is not connected, contradiction. This proves (i).
For (ii), suppose, to the contrary, that there exists $F \not= F'$ from $\mathfrak{C}$ with $F \cap S= \emptyset$
and $F' \cap S = \emptyset$. By (i), there exist $x \in F\cap V(C)$ and $y \in F' \cap V(D)$, but no $x$,$y$-path in $G$ avoiding $S$
and, hence, no $x$,$y$-path in $G[F \cup F']$, contradiction.
\hspace*{\fill}$\Box$

In particular, every graph with a Kempe coloring of order $k$ must be $(k-1)$-connected.
We repeat the following Lemma (and its proof) from~\cite{Kriesell2001}. 

\begin{lemma}~\cite{Kriesell2001}
  \label{L2}
  Suppose that $H$ is a graph such that $L(H)$ is $k$-connected. Then for all distinct vertices $a, b$ of degree at least $k$, there
  exist $k$ edge-disjoint $a$,$b$-paths in $H$.
\end{lemma}

{\bf Proof.}
If there were no such paths, then, by Menger's Theorem (cf.~\cite{Diestel2017}),
there exists an $a,b$-cut $S$ in $H$ with less than $k$ edges.
We may assume that $S$ is a minimal $a,b$-cut, implying that $H-S$ has exactly two components $C,D$,
where $a \in V(C)$ and $b \in V(D)$. Since both $a,b$ have degree at least $k$,
there exists an edge $e \in E_H(a) \setminus S$, that is, $e \in E(C)$, and at least one edge $f \in E(D)$.
But then $S$ separates $e$ from $f$ in $H$ and, thus, $e$ from $f$ in $L(H)$, contradicting the assumption that $L(H)$ is $k$-connected.
\hspace*{\fill}$\Box$

We are now ready to prove our main result. At some places, we will contract some subgraph $X$ of $H$ to a single vertex.
In order to make object references easier, we choose a graph model where
the edge set of the resulting graph actually {\em equals} $E(H) \setminus E(X)$,
not just ``corresponds  to $E(H) \setminus E(X)$'' in whatever way.

{\bf Proof of Theorem~\ref{T2}.}

We proceed by induction on $|E(H)|$. 
Let $\mathfrak{C}$ be a partition of $E(H)$ into matchings such that the union of any two of them is connected,
set $k:=|\mathfrak{C}|$, and let $T$ be a transversal of $\mathfrak{C}$.
We have to show that there exists a set of $k$ connected, pairwise disjoint, pairwise incident edge sets traversed by $T$ in $H$.
Observe that the statement is easy to prove whenever $k \leq 2$. Hence, we may assume $k \geq 3$. In particular, $|E(H)|\geq 3$. 
By Lemma~\ref{T3} we may assume that $H$ is simple.
Observe that $\mathfrak{C}$ is a Kempe coloring of $L(H)$.

As $\mathfrak{C}$ is a partition of $E(H)$ into $k$ matchings, all vertices in $H$ have degree at most $k$.
Suppose first that there is a vertex $v$ of degree $k$ in $H$.
Then $U:=E_H(v)$ induces a clique of order $k$ in $L(H)$.
If there is a set of $k$ disjoint $U$,$T$-paths in $L(H)$, then their vertex sets form a clique minor in $L(H)$
and, at the same time, a set of $k$ connected, pairwise disjoint, pairwise incident edge sets in $H$, traversed by $T$.
So we may assume that there are no $k$ disjoint $U$,$T$-paths in $L(H)$.
By Menger's Theorem (cf.~\cite{Diestel2017}), there exists a vertex set $S$ in $L(H)$ with $|S|<k$ separating $U$ from $T$,
and we may take a smallest such set, implying that $L(H)-S$ has only two components $C,D$
(since for every vertex $v \in S$, $N_{L(H)}(v)$ contains a vertex from each component of $L(H)-S$ but is, at the same
time, the union of at most two cliques of $L(H)$). 
We may assume
that, say, $C$ contains at least one vertex from $U$ as $|U|>|S|$. But then $C$ contains no vertex from $T$,
and $D$ contains at least one vertex from $T$ and no vertex from $U$. 
Back in $H$,  the set $S$ is a minimal cut in $H$
and $H-S$ has exactly two components, $C',D'$, where $E(C')=V(C)$ and $E(D')=V(D)$. 
As $E(C')$ contains an edge from $U$, $v \in V(C')$ follows.

From Lemma~\ref{L1} we know that $S$ consists of $k-1$ objects, all coming from distinct members of $\mathfrak{C}$.
Let $F$ be the unique member from $\mathfrak{C}$ with $F \cap S = \emptyset$.
Let $H'$ be the graph obtained from $H$ by contracting $D'$ to a single vertex $w$. 
Since $F \in \mathfrak{C}$ contains an edge from $E(C')$ by Lemma~\ref{L1} and all other classes contain an edge from $S$,
we obtain a partition of $E(H')$ into matchings such that the union of any two of them is connected,
by deleting all edges of $E(D')$ from their sets in $\mathfrak{C}$.
Thus, we also have a Kempe coloring of $L(H')$, so that, by Lemma~\ref{L1}, $L(H')$ is $(k-1)$-connected.
Since $v$ has degree at least $k$ and $w$ has degree $k-1$ in $H'$, there exist $k-1$ edge-disjoint $v$,$w$-paths in $H'$ by Lemma~\ref{L2}.
For $e \in S$, let $P_e$ be the path among these containing $e$.

Now let $H''$ be the graph obtained from $H$ by contracting $C'$ to a single vertex (recall that $E(C')\neq \emptyset$). As above, by deleting all edges of $E(C')$ from their
sets in $\mathfrak{C}$, we obtain a partition of $E(H'')$ into matchings such that the union of any two of them is connected.
Moreover, $T$ remains a transversal of the modified partition.
Since $|E(H'')|<|E(H)|$, we may apply induction to $H''$ and find a set $\mathfrak{K}$ of connected, pairwise disjoint, pairwise incident edge sets
traversed by $T$ in $H''$.
Setting $A':=A \cup \bigcup \{E(P_e):\, e \in A,\,e \in S\}$ for $A \in \mathfrak{K}$, one readily checks that $\mathfrak{K}':=\{A':\,A \in \mathfrak{K}\}$
proves the statement of the theorem for $H$.

Therefore, we may assume from now on that the maximum degree $\Delta$ of $H$ is at most $k-1$. Let $\delta$ denote the minimum degree of $H$.
Every edge $xy$ is incident with at least one edge from each of the $k-1$ members of $\mathfrak{C}$ not containing $xy$,
so that $d_H(x)+d_H(y) \geq k+1$. Consequently, $\delta \geq k+1 - \Delta$ and $\Delta \geq (k+1)/2$.
For distinct $A,B$ from $\mathfrak{C}$, consider the subgraph $H(A,B)$ formed by all edges of $A \cup B$.
It is either a path or a cycle, and we say that $H(A,B)$ {\em ends} in a vertex $v$ if $v$ has degree $1$ in $H(A,B)$,
or alternatively, if $v$ is covered by exactly one edge of $A,B$.
Now, if $v$ has degree $d$ in $H$, then it is covered by exactly $d$ of the $k$ matchings from $\mathfrak{C}$,
so that exactly $d \cdot (k-d)$ of the subgraphs $H(A,B)$ end in $v$.
Observe that $k-\Delta<\delta \leq d \leq \Delta$ and
consider the quadratic function $f$ defined by $f(d):=d \cdot (k-d)-\Delta \cdot (k-\Delta)$ with zeroes $\Delta$ and $k-\Delta$.
We get $f(d) \geq 0$, that is, $d \cdot (k-d) \geq \Delta \cdot (k-\Delta)$, for $d \in (k-\Delta,\Delta]$ with equality only if $d=\Delta$.
Consequently, in each vertex $v$, for at least $\Delta \cdot (k-\Delta)>0$ pairs $A,B\in \mathfrak{C}$, the graph $H(A,B)$ ends in $v$. As there are only ${k \choose 2}$
many subgraphs $H(A,B)$ and as each of them ends in two or zero vertices, we get $|V(H)| \cdot \Delta \cdot (k-\Delta) \leq k \cdot (k-1)$.
Moreover, $|V(H)| \geq \Delta+1$ since $H$ is simple, so
\begin{eqnarray}
\label{I1}
(\Delta+1) \cdot \Delta \cdot (k-\Delta) \leq k \cdot (k-1).
\end{eqnarray}
Consider the cubic function $g$ defined by $g(\Delta):=(\Delta+1) \cdot \Delta \cdot (k-\Delta) - k \cdot (k-1)$.
It has zeros $k-1$ and $\pm \sqrt{k}$, so that it is positive for $\Delta \in (\sqrt{k},k-1)$, that is, $(\Delta+1) \cdot \Delta \cdot (k-\Delta) > k \cdot (k-1)$
for all $\Delta \in [(k+1)/2,k-1)$ as $k\geq 3$. 
Since $\Delta \in [(k+1)/2,k-1]$ and (\ref{I1}) holds, this necessarily implies $\Delta=k-1$ and equality in (\ref{I1}). Backtracking through the
arguments leading to (\ref{I1}) yields, subsequently: $|V(H)|=\Delta+1$; in each vertex $v$, exactly $\Delta \cdot (k-\Delta)$ of the $H(A,B)$ end;
and, finally, each vertex of $H$ has degree $d=\Delta$.

It follows that $H$ is the simple complete graph on $\Delta+1=k \geq 3$ vertices, and we obtain the statement of Theorem~\ref{T2} for $H$ from Lemma~\ref{T4}.
\hspace*{\fill}$\Box$ 

\bigskip

\centerline{*}

\medskip

\centerline{*}

{\bf Authors' addresses.}

Matthias Kriesell and Samuel Mohr \\
Institut f\"ur Mathematik der Technischen Universit\"at Ilmenau \\
Weimarer Stra{\ss}e 25 \\
98693 Ilmenau \\
Germany

\end{document}